\documentclass[12pt]{article}
\usepackage[intlimits]{amsmath}
\usepackage{amssymb,mathrsfs}
\usepackage{indentfirst}

\renewcommand{\le}{\leqslant}                                                    
\renewcommand{\ge}{\geqslant}

\title{A new randomized algorithm for the Erd\H{o}s--Hajnal problem}

\author{D.D. Cherkashin\footnote{Saint-Petersburg State University, Department of Mathematics and Mechanics; e-mail: matelk@mail.ru.}}

\date{}

\begin{document}

\maketitle

\begin{abstract}

In 1961 Erd\H{o}s and Hajnal introduced the quantity $m(n)$ as the minimum number of edges in an $n$-uniform hypergraph with chromatic number at least 3. 
The best known lower and upper bounds for $ m(n) $ are $ c_1 \sqrt{\frac{n}{\ln n}} 2^n$ and $c_2 n^2 2^n$ respectively. The lower bound is due 
to Radhakrishnan and Srinivasan (see \cite{RS}). 
A natural generalization for $ m(n) $ is the quantity $ m(n,r) $, which is the minimum number of edges in an $n$-uniform hypergraph with chromatic number at 
least $r+1$. In this work, we present a new randomized algorithm yielding a bound $ m(n,r) \ge c n^{\frac{r-1}{r}} r^{n-1} $, 
which improves upon all the previous bounds in a wide range of the parameters $ n, r $. Moreover, for $ r = 2 $, we get exactly the same bound as 
in the work \cite{RS} of Radhakrishnan and Srinivasan, and our proof is simpler. 

\end{abstract}


\section{Introduction}

In 1961 Erd\H{o}s and Hajnal (see \cite{EH}) introduced the quantity $ m(n) $, which is the minimum number of edges in an $ n $-uniform hypergraph 
$ H $ with chromatic number greater than 2. The problem of finding the right order of $ m(n) $ is one of the most important in extremal combinatorics, and 
substantial work has been done in this direction (see, e.g., \cite{Kost1}, \cite{Rai1}). In particular, the best known lower bound for the value $ m(n) $
is due to Radhakrishnan and Srinivasan (see \cite{RS}) who proposed in 2000 an ingenious randomized algorithm resulting in the following theorem. 

\vskip+0.2cm

\noindent {\bf Theorem 1.} {\it For any $ q < \sqrt{2} $, there exists $ n_0 $ such that for every $ n \ge n_0 $, the inequality holds
$$
m(n) \ge q \left(\frac{n}{\ln n}\right)^{1/2} 2^{n-1}.
\eqno{(1)}
$$} 

\vskip+0.2cm

An important generalization of the quantity $ m(n) $ is the quantity $ m(n,r) $, which is the minimum number of edges in an $ n $-uniform hypergraph 
$ H $ with chromatic number greater than $ r $ (see \cite{HS}). For $ m(n,r) $, several different lower bounds are found, each of them being currently 
the strongest for some range of the parameters $ n, r $. Let us briefly describe the whole picture. First, we mention the work \cite{Sh1} by Shabanov 
who just noticed that the randomized algorithm by Radhakrishnan and Srinivasan can be generalized to give the bound 
$$
m(n,r) \ge (\sqrt{3}-1)\left(\frac n{\ln n}\right)^{1/2}r^{n-1}
\eqno{(2)}
$$
for any $ n \ge 2 $, $ r \ge 2 $. Second, there was an interesting and rather involved algorithm by Kostochka (see \cite{Kost2}) proving that for 
$ r = O\left(\sqrt{\ln \ln n}\right) $, one has
$$
m(n,r) \ge e^{-4r^2} \left(\frac n{\ln n}\right)^{a/(a+1)}r^n,
\eqno{(3)}
$$
where $a=\left\lfloor \log_2 r\right\rfloor$. Finally, there were some improvements by Pluh\'ar and Shabanov (see \cite{Pluh}, \cite{Sh2}): 
$$
m(n,r)\ge \frac 14\; n^{1/2} r^{n-1}, \,\,\, n \ge 2, \,\,\, r \ge 3,
\eqno{(4)}
$$
and
$$
m(n,r)\ge \left(\pi^{\frac 1r}\;e^{-\frac 1{12(n-1)}}\right)\frac 1{e\sqrt{2\pi}}\; (n-1)^{\frac 12-\frac 1{2r}}\;r^{n+2/r}.
\eqno{(5)}
$$
It is worth noting here that sometimes we write in the bounds $ r^n $ and sometimes we write $ r^{n-1} $. The difference becomes substantial if 
$ r $ is a growing function of $ n $. 

In this paper we propose a new randomized algorithm, and using it we prove the following theorem. 

\vskip+0.2cm

\noindent {\bf Theorem 2.} {\it Let $ r = o\left(\frac{\ln n}{\ln\ln n}\right) $. Let 
$$
q < \frac{(r!)^{\frac{1}{r}}}{2^{\frac{1}{r}} \left(1-\frac{1}{r}\right)^{\frac{r-1}{r}} r^{\frac{r-2}{r}}} = : q'.
$$ 	
Then there exists $ n_0 $ such that for every $ n \ge n_0 $, the inequality holds
$$
m(n,r) \ge q \left(\frac{n}{\ln n}\right)^{\frac{r-1}{r}} r^{n-1}.
\eqno{(6)}
$$} 

\vskip+0.2cm

On the one hand, estimate (6) improves substantially upon all the previous estimates (2)--(5), provided $ r $ is not too large. On 
the other hand, if we put $ r=2 $ in Theorem 2, then we get {\it exactly} the result of Theorem 1. As we will see from the proof in the next section, 
our proof is even shorter than that of Theorem 1.  

Note that the condition $ r = o\left(\frac{\ln n}{\ln\ln n}\right) $ is sufficient, but not necessary. It is just taken for simplicity of the proof. 

To complete this section we present one more result, which will be obtained automatically from the proof of Theorem 2. 

\vskip+0.2cm

\noindent {\bf Theorem 3.} {\it Let $ r = o\left(\frac{\ln n}{\ln\ln n}\right) $. Then there exists a constant $ c > 0 $ such that whenever an 
$ n $-uniform hypergraph has maximum edge degree 
$$
D < c \left(\frac{n}{\ln n}\right)^{\frac{r-1}{r}} r^{n-1},
$$
its chromatic number is not greater than $ r $.}

\vskip+0.2cm

\section{Proofs of new results}

\subsection{An auxiliary result}

In this section, we obtain a result, which is suggested by Pluh\'ar's work \cite{Pluh}. First, we briefly describe Pluh\'ar's assertion. Let $ H = (V,E) $ be 
a hypergraph (not necessarily uniform). Let $ \sigma $ denote an ordering of the vertices $ V $. A family of edges $ A_1, \dots, A_r \in E $ is called 
{\it ordered $ r $-chain} in $ \sigma $, if the following conditions hold:

\begin{enumerate}

\item for any $ i \in \{1, \dots, r-1\} $, we have $ |A_i \cap A_{i+1}| = 1 $;

\item for any $ i, j $ such that $ |i-j| > 1 $, we have $ A_i \cap A_{i+1} = \emptyset $; 

\item for any $ i \in \{1, \dots, r-1\} $ and any $ v \in A_i $, $ u \in A_{i+1} $, we have $ \sigma(v) \le \sigma(u) $.

\end{enumerate}

Proposition 1 below was proved by Pluh\'ar. 

\vskip+0.2cm

\noindent {\bf Proposition 1.} {\it The chromatic number of a hypergraph $ H = (V,E) $ does not exceed $ r $ if and only if there exists an ordering $ \sigma $ of 
$ V $ such that there are no ordered $ r $-chains in $ H $.} 

\vskip+0.2cm

Here we introduce the notion of a {\it strong ordered $ r $-chain}. To this end we assume that there is a map $ f $ from the set of edges $ E $ of a hypergraph 
to the set $ \{1, \dots, r\} $, i.e., to each edge a label is assigned. A family of edges $ A_1, \dots, A_r \in E $ is called 
{\it strong ordered $ r $-chain} in $ \sigma $, if it is subject to conditions 1--3 and, moreover, for any $ i $, we have $ f(A_i) = i $. Note that both notions 
of $r$-chains have sense even in the cases of hypergraphs with multiple edges and with edges of length 1. The only subtlety here is as follows: if we have 
coinciding edges of length 1, then condition 2 is not completely correct; for example, a chain $ \{1,2\}, \{2\}, \{2\}, \{2,3\} $ is admissible. 

\vskip+0.2cm

\noindent {\bf Proposition 2.} {\it Consider a hypergraph $ H = (V,E) $ and a map $ f $ from $ E $ to $ \{1, \dots, r\} $. Then the following statements are equivalent:
(i) there is a coloring of $ V $ in $ r $ colors such that no edge $ A $ in $ E $ consists only of vertices of color $ f(A) $;
(ii) there is an order of elements of $ V $ without strong ordered $r$-chains.} 

\vskip+0.2cm

\paragraph{Proof of Proposition 2.} {\it (ii) implies (i).} Fix an ordering $\sigma$ of $V$. Make a greedy coloring. First, color all the vertices in color 1. 
Then consider the vertices in order $\sigma$. If the $i$th vertex is the first vertex of some edge, then recolor it in the color with 
the minimal number such that there is no conflict with edges that contain only vertices from the set $1, \dots, i$.
If this is impossible, then recolor the $i$th vertex with color $r$. 

Clearly if the final coloring is not good, then any bad edge has color $r$.
Consider an arbitrary bad edge $A_r$. By definition $f(A_r) = r$.
The first vertex of $A_r$ could not be colored $r-1$. Hence it is the last vertex of an edge $A_{r-1}$ all of whose 
vertices except for one are colored $r-1$; moreover, $f(A_{r-1}) = r-1$. Going this way we get a strong ordered $r$-chain.

{\it (i) implies (ii).} Assume we have a good coloring. Then we order the vertices by growth of color. If there is a strong ordered $r$-chain, then it is 
easy to prove by induction that the $i$th edge in it has a vertex whose color is greater than $ i $, which is impossible for the last edge. It is a 
contradiction. 

\vskip+0.2cm

\subsection{Proof of Theorem 2}

Fix a hypergraph $ H = (V,E) $ with 
$$
|E| \le q \left(\frac{n}{\ln n}\right)^{\frac{r-1}{r}} r^{n-1}.
$$
We want to show that its chromatic number is not greater than $ r $. As it is usual in this type of proofs, we make a random coloring. The procedure 
consists of two steps. At the first step, we color every vertex (independently) into color $ i \in \{1, \dots, r\} $ 
with probability $ \frac{1-p}{r} $ for some $ p $. So with 
probability $ p $ a vertex is colorless. Suppose there are no monochromatic, almost monochromatic, and fully colorless edges (an edge is almost 
monochromatic, if exactly one of its vertices is colorless and all the other vertices have same color). 
At the second step, we take the set $ W \subseteq V $ of colorless vertices. Then for every $ A \in E $, we 
take $ B = A \cap W $: $ |B| \not \in \{0, 1, |A|\} $. We get a hypergraph $ H' = (W,F) $. If an edge $ B \in F $ is such that $ A \setminus B $ is not 
monochromatic after the first step, we may forget about it, since it is already colored in a right way. Thus, we consider only those $ B \in F $, 
for which $ A \setminus B $ is completely colored with some color $ i \in \{1, \dots, r\} $. Let $ f(B) = i $ and take a random ordering of $ W $. 
If there are no strong ordered $ r $-chains, then by Proposition 2 we have a good coloring. 

We must prove that the probability of a failure is less than 1. Put
$$
p = \frac{(r-1)\ln \left(r^{\frac{r}{r-1}} n \ln^2 n\right)}{rn}.
$$ 

The probability that a monochromatic edge appears after the first step is bounded from above by 
$$
|E| \cdot r \cdot \left(\frac{1-p}{r}\right)^n \le |E| r^{1-n} n^{\frac{1-r}{r}} r^{-1} \left(\ln^2 n\right)^{\frac{1-r}{r}} \le 
q \left(\frac{n}{\ln n}\right)^{\frac{r-1}{r}} r^{n-1} r^{1-n} n^{\frac{1-r}{r}} r^{-1} \left(\ln^2 n\right)^{\frac{1-r}{r}} \to 0, \,\,\, 
n \to \infty. 
$$
The probability that an almost monochromatic edge appears after the first step is bounded from above by 
$$
|E| \cdot rnp \cdot \left(\frac{1-p}{r}\right)^{n-1} \le 
(1+o(1)) q \left(\frac{n}{\ln n}\right)^{\frac{r-1}{r}} r^{n-1} r^{2-n} (\ln n) 
n^{\frac{1-r}{r}} r^{-1} \left(\ln^2 n\right)^{\frac{1-r}{r}} \sim q (\ln n)^{1-3\frac{r-1}{r}} \to 0. 
$$
The probability that a fully colorless edge appears after the first step does not exceed $ |E| p^n $, which is tiny. 

Fix an $r$-chain in $ H' $. Let $ a_1, \dots, a_r \ge 2 $ be the cardinalities of its edges. The probability that this chain is ordered equals
$$
2 \frac{(a_1-1)!(a_r-1)!\prod\limits_{i=2}^{r-1} (a_i - 2)!}{\left(\sum a_i - r + 1\right)!} =: M(a_1, \dots, a_r).
$$

The probability that our chain is strong equals
$$
p^{r-1} p^{a_1-1} \left(\frac{1-p}{r}\right)^{n-a_1} C_{n-1}^{a_1-1} 
\left(\prod_{i=2}^{r-1} p^{a_i-2} \left(\frac{1-p}{r}\right)^{n-a_i} C_{n-2}^{a_i-2}\right) p^{a_r-1}\left(\frac{1-p}{r}\right)^{n-a_r} C_{n-1}^{a_r-1}
=: N(a_1, \dots, a_r).
$$

The number of $ r $-chains is less than or equal to $ |E|^r/r! $. Thus the final probability does not exceed the value
$$
\frac{|E|^r}{r!} \sum_{a_1, \dots, a_r} M(a_1, \dots, a_r) N(a_1, \dots, a_r) \le 
\frac{2|E|^r}{r!} \sum_{a_1, \dots, a_r} \frac{p^{\sum a_i - r + 1} n^{\sum a_i - 2r +2} \left(\frac{1-p}{r}\right)^{nr - \sum a_i}}{\left(\sum a_i - r + 1\right)!}.
$$
Fix $ t \in {\mathbb N} $. The number of collections of $ a_i $'s with $ \sum a_i - 2r + 2 = t $ is less than $ t^{r-1} $. Hence, we may write the estimate
$$
\frac{2|E|^r}{r!} \sum_{t=0}^{\infty} t^{r-1} 
\frac{p^{t+r-1} n^{t} \left(\frac{1-p}{r}\right)^{nr - t - 2r + 2}}{\left(t + r - 1\right)!} \le 
\frac{2|E|^r}{r!} p^{r-1} \left(\frac{1-p}{r}\right)^{nr - 2r + 2} \sum_{t=0}^{\infty}  
\frac{p^{t} n^{t} \left(\frac{1-p}{r}\right)^{-t}}{t!} \le
$$
$$
\le \frac{q}{q'} (1+o(1)) \left(1-\frac{1}{r}\right)^{1-r} r^{2-r} \left(\frac{n}{\ln n}\right)^{r-1} r^{nr-r} 
\left(1-\frac{1}{r}\right)^{r-1} n^{1-r} (\ln n)^{r-1} e^{-pnr} r^{2r-2-nr} e^{\frac{pnr}{1-p}} = \frac{q}{q'} (1+o(1)).
$$

Theorem 2 is proved. 

\subsection{Proof of Theorem 3}

Let us proceed with the same random coloring as in the proof of Theorem 2. We have the two following types of events. On the one hand, at the first step of 
our algorithm, we can get a monochromatic, an almost monochromatic, or a fully colorless edge. Let $ P_1 $ be the probability of any such event (corresponding
to a fixed edge). On the other hand, at the second step, we can get a strong ordered $r$-chain. Let $ P_2 $
be the probability of any such event (corresponding to a fixed $r$-chain). The values of $ P_1, P_2 $ are given in the previous section. 

We have to prove that none of these events occurs with positive probability. It is 
quite standard to use here Lov\'asz Local Lemma (LLL). Clearly any edge intersects at most $ D $ other edges and $ D^r $ different $ r $-chains. Also any $ r $-chain 
intersects at most $ rD $ edges and $ rD^r $ other $ r $-chains. Thus, to utilize LLL we just need to find such $ x $ and $ y $ that simultaneously satisfy the 
inequalities
$$
P_1 \le x (1-x)^D (1-y)^{D^r}, \,\,\,\,\, P_2 \le y (1-x)^{rD} (1-y)^{rD^r}.
$$

A tedious and standard calculation that we omit here shows that taking $ x = \frac{a}{D} $ with an arbitrary positive $ a $ and $ y \asymp \frac{1}{rD^r} $ one 
eventually gets the result of Theorem 3.

\end{document}